# Lucas sequences whose $n$th term is a square or an almost square

A. Bremner, N. Tzanakis

6 October 2006

**Abstract**

Let $P$ and $Q$ be non-zero integers. The Lucas sequence $\{U_n(P,Q)\}$ is defined by $U_0 = 0$, $U_1 = 1$, $U_n = PU_{n-1} - QU_{n-2}$, $(n \geq 2)$. Historically, there has been much interest in when the terms of such sequences are perfect squares (or higher powers). Here, we summarize results on this problem, and investigate for fixed $k$ solutions of $U_n(P,Q) = k\Box$, $(P,Q) = 1$. We show finiteness of the number of solutions, and under certain hypotheses on $n$, describe explicit methods for finding solutions. These involve solving finitely many Thue-Mahler equations. As an illustration of the methods, we find all solutions to $U_n(P,Q) = k\Box$ where $k = \pm 1, \pm 2$, and $n$ is a power of 2.

## 1 Introduction

Let $P$ and $Q$ be non-zero integers. The Lucas sequence $\{U_n(P,Q)\}$ is defined by

$$U_0 = 0, \quad U_1 = 1, \quad U_n = PU_{n-1} - QU_{n-2} \quad (n \geq 2). \tag{1}$$

Historically, there has been much interest in when the terms of such sequences are perfect squares (or higher powers), and we summarize here the numerous and diverse results. Ljunggren [9] shows that for $(P,Q) = (2,-1)$ and $n \geq 2$, then $U_n$ is a perfect square precisely for $U_7 = 13^2$, and $U_n = 2\Box$ precisely for $U_2 = 2$. The sequence $\{U_n(1,-1)\}$ is the familiar Fibonacci sequence, and Cohn [5] proved in 1964 that the only perfect square greater than 1 in this sequence is $U_{12} = 12^2$. Ribenboim and McDaniel [17] show that for $P$ even and $Q \equiv 1 \bmod 4$, then $U_n(P,Q) = \Box$ imposes necessary conditions on the prime factorization of $n$. Earlier, in [16], the same authors show with only elementary methods that when $P$ and $Q$ are *odd*, and $P^2 - 4Q > 0$, then $U_n = \Box$ only for $n = 0,1,2,3,6$ or $12$; and that there are at most two indices greater than 1 for which $U_n$ can be square. They characterize fully the instances when $U_n = \square$, for $n = 2, 3, 6$. Bremner & Tzanakis [2] extend these results by determining all Lucas sequences $\{U_n(P, Q)\}$ with $U_{12} = \square$, subject only to the restriction that $\gcd(P, Q) = 1$ (it turns out that the Fibonacci sequence provides the only example). Under the same hypothesis, all Lucas sequences with $\{U_n(P, Q)\}$ with $U_9 = \square$ are determined. In a later paper, the same authors [3] show that for $n = 2, ..., 7$ then $U_n(P, Q)$ is square for infinitely many coprime $P, Q$ and determine all sequences $\{U_n(P, Q)\}$ with $U_n(P, Q) = \square$, $n = 8, 10, 11$.

We discuss in this paper the more general problem of finding all integers $n$, $P, Q$, for which $U_n(P, Q) = k\square$ for a given integer $k$. Results of Pethő [14], Shorey and Stewart [19], and Shorey and Tijdeman [20], show the finiteness of solutions of $U_n(P, Q) = k\square$ for fixed $P, Q$. Actually, these results are much more general as their scope is not only squares and Lucas sequences, but perfect powers in general, in a much larger class of second order recurrence sequences. The results are effective, but not however explicit. Yabuta [23] gives some special results concerning squares in the Lucas sequences (effective, but not explicit), proved by elementary means. When $k$ is restricted to the possibilities $1, 2, 3, 6$ several authors have given specific results as follows. Mignotte and Pethő [10] show that if $P \geq 3$, then $U_n(P, 1) = \square$ for $n \geq 3$ has exactly the solutions $U_4(338, 1) = 6214^2$, $U_6(3, 1) = 12^2$; and further, again for $P \geq 3$, that $U_n(P, 1) = k\square$, $k = 2, 3, 6$, has no solutions for $n \geq 4$. Nakamula and Pethő [12] show that for $P \geq 1$ and $n \geq 4$, then $U_n(P, -1) = m^2$ implies $(P, n, m) = (1, 12, 12)$, $(2, 7, 13)$; that $U_n(P, -1) = 2m^2$ implies $(P, n, m) = (1, 6, 2)$, $(4, 4, 6)$; that $U_n(P, -1) = 3m^2$ implies $(P, n, m) = (1, 4, 1)$, $(2, 4, 2)$, $(24, 4, 68)$; and $U_n(P, -1) = 6\square$ has no solutions.

If $n_0 \geq 8$ is fixed, then for fixed $k$ it is shown here that there are at most finitely many coprime $P, Q$ for which $U_{n_0}(P, Q) = k\square$ (in fact, for $n_0 \geq 9$ this follows trivially from either Faltings' Theorem or from Darmon and Granville [7]). More generally, we consider the problem of *explicitly* determining all integers $n = mn_0$ and $P, Q$ such that $U_n(P, Q) = k\square$, under the hypothesis that the prime divisors of $m$ belong to a prescribed finite set of primes. To put this problem in context and give some general motivation, we make the following remarks about the equation $u_n(P, Q) =$"almost" r-th power, where $u_n$ is a second order recurrence sequence defined by $u_n = Pu_{n-1} - Qu_{n-2}$. The associated problems, of types (i)-(vi), are summarised in the following table, in which the indication $c$ denotes known, the indication $x$ denotes unknown. The final two columns indicate bibliographic references for the appropriate problem type.

| Type | $(P, Q)$ | $n$ | $r$ | References | Comments |
|---|---|---|---|---|---|
| i | $c$ | $x$ | $c$ | [5], [9] | |
| ii | $c$ | $x$ | $x$ | [4] | an explicit result, highly technical |
| | | | | [14], [19], [20] | finiteness results |
| iii | $x$ | $c$ | $c$ | [2], [3] | |
| iv | $x$ | $x$ | $c$ | [10], [12] | $P$ unknown, $Q$ prescribed as $\pm 1$ |
| | | | | [16], [17] | $P, Q$ restricted, elementary methods |
| | | | | [23] | finiteness results |
| v | $x$ | $c$ | $x$ | | |
| vi | $x$ | $x$ | $x$ | | |

Problems of type (i) in general lead to explicit Thue-Mahler equations, with associated intricacy of computation (the above papers are exceptions, [5] using only elementary methods, and [9] avoiding Thue-Mahler equations). See Section 4 of this paper for examples of type (i) with $r = 2$. Problem (ii) is very difficult. A solved example of this type may be found in Bugeaud, Mignotte and Siksek [4]. Problem (iii) leads to solving hyperelliptic or superelliptic equations. Bremner and Tzanakis [2], [3], solve problems of this type. Problems (iv), (v), (vi) are of essentially strictly increasing difficulty. Problem (iv) is already very difficult, though may be treated by elementary methods under certain restrictions on $P, Q$ such as in [16], [17]. The problem referred to above that we study in this paper (viz. Problem 1 of Section 2) is a restriction of Problem (iv). Its solution reduces to two separate steps, which may be succinctly described as follows: first, that of finding the finitely many $P, Q$ satisfying $U_{n_0}(P, Q) = k_0\Box$, where $k_0$ runs through a finite set of square-free integers whose prime divisors are those of $m$ and/or $k$; and second, for each such pair $P_0, Q_0$, finding all integers $n$ with $n_0|n$ and $U_n(P_0, Q_0) = k\Box$ (in fact to find all such $n$ at this stage, the requirement that $n_0|n$ can in principle be relaxed). The first step demands finding explicitly all rational points on certain curves of genus greater than 1. Although this task has in general not yet been proved to be effective, in any particular example of "reasonable" nature, there is expectation that an explicit determination of all points may be obtained. For example, the papers [2], [3] perform this computation in the case $k = 1$, $8 \leq n \leq 12$; see also the relevant bibliography therein. The second step is shown to reduce to solving finitely many Thue-Mahler equations, which, at least in principle, are explicitly solvable; see Tzanakis and de Weger [21].

The layout of the paper is as follows. Sections 2, 3 show that Problem 1 has an effective solution. In Section 2, finiteness of the number of solutions $(P, Q)$ to $U_{n_0}(P,Q) = k\square$, $n_0 \geq 8$, is established (of course, explicitly finding these solutions may be very difficult). The main result is Theorem 2.1. In Section 3, we show that for fixed $P, Q, k$, finding $n$ such that $U_n(P,Q) = k\square$ leads to solving finitely many quartic Thue-Mahler equations, which in principle may be effectively solved. Section 4 illustrates the methods by explicitly finding all $P, Q, n$ with $U_n(P,Q) = k\square$, where $k = \pm 1, \pm 2$ and $n$ is a power of 2. The main result is Theorem 4.2. Its proof requires Theorem 4.1 (which solves problems of type (iii)), with proof of this Theorem given in Section 5. It also requires solving two Thue-Mahler equations, namely $u^4 - 17v^4 = \pm 2^{z+2}$ and $u^4 - 84v^4 = 17^z$, which are solved in Sections 6 and 7, respectively.

## 2 Preliminaries

We consider the Lucas sequence with parameters $P, Q$, defined by

$$U_0(P,Q) = 0\,, \quad U_1(P,Q) = 1\,, \quad U_n = PU_{n-1} - QU_{n-2} \quad \text{for } n \geq 2$$

and the associated Lucas sequence

$$V_0(P,Q) = 2\,, \quad V_1(P,Q) = P\,, \quad V_n = PV_{n-1} - QV_{n-2} \quad \text{for } n \geq 2\,.$$

The characteristic polynomial $t^2 - Pt + Q$ has discriminant $D = P^2 - 4Q$. Throughout this paper we assume that

$$PQ \neq 0\,, \quad \gcd(P,Q) = 1 \quad \text{and } D \neq 0. \tag{2}$$

We will often use various properties of the Lucas sequences, which can be found, for example, in the first two sections of [15].

For any non-zero integer $x$ we define

$$\mathcal{P}(x) = \begin{cases} \{1\} & \text{if } |x| = 1 \\ \text{Set of primes dividing } x & \text{if } |x| > 1 \end{cases}.$$

**Problem 1.** *Let $k \neq 0$ and $n_0 > 1$ be fixed integers. Let $\mathcal{S}$ be a fixed non-empty finite set of primes. Find all integral triads $(m, P, Q)$, under the constraint (2), for which*

$$\mathcal{P}(m) \subseteq \mathcal{S}\,, \quad U_{n_0 m}(P,Q) = k\square\,.$$

Until the end of this section we will often write, for example, $U_n$ instead of $U_n(P,Q)$, in order to avoid overloaded notation. Put $n_0 m = n$ and, of course, $U_{n_0}|U_n$; also, put $d = \gcd(U_{n_0}, U_n/U_{n_0})$. By results of Lehmer [8] (see also (2.2) of [15]), $d|m$,[1] hence $\mathcal{P}(d) \subseteq \mathcal{P}(m) \subseteq \mathcal{S}$. On the other hand, the relation $U_n = k\Box$ can be written

$$\frac{U_{n_0}}{d} \cdot \frac{U_n/U_{n_0}}{d} = k\Box\,,$$

where the factors on the left-hand side are relatively prime; hence

$$U_{n_0}/d = \pm k_0\Box\,, \tag{3}$$

where $k_0$ is a positive integer with $\mathcal{P}(k_0) \subseteq \mathcal{P}(k)$. Equation (3) therefore reduces to

$$U_{n_0}(P,Q) = \pm p_1^{i_1}\cdots p_r^{i_r}\Box\,, \quad i_1,\ldots,i_r \in \{0,1\}\,, \tag{4}$$

where

$$\{p_1,\ldots,p_r\} \subseteq \mathcal{P}(k) \cup \mathcal{P}(d) \subseteq \mathcal{P}(k) \cup \mathcal{S}\,.$$

In general, if $n_0 \geq 8$, then any equation (4) has only finitely many solutions in relatively prime integers $P, Q$. Indeed, we have the following result.

**Theorem 2.1.** *Let $k, n_0$ be fixed non-zero integers with $n_0 \geq 8$. Then $U_{n_0}(P,Q) = k\Box$ can hold for at most finitely many relatively prime integers $P, Q$.*

*Proof.* For $n_0$ odd, we have $U_{n_0}(P,Q) = P^{n_0-1}U_{n_0}(1, \frac{Q}{P^2})$, so that $U_{n_0}(P,Q) = k\Box$ implies that the curve $U_{n_0}(1,x) = ky^2$ contains a rational point with $x = \frac{Q}{P^2}$. For $n_0$ even, we have by induction that $U_{n_0}(P,Q) = PU'_{n_0}(P,Q)$, where $U'_{n_0}(P,Q) \in \mathbb{Z}[P,Q]$, with $U'_{n_0}(P,Q)$ homogeneous in $P^2, Q$ of degree $\frac{n_0}{2} - 1$. Also by induction, $U'_{n_0}(P,Q) \equiv \frac{n_0}{2}(-Q)^{\frac{n_0}{2}-1} \bmod P$ (indeed, $\bmod P^2$), so that the greatest common divisor of $P$ and $U'_{n_0}(P,Q)$ divides $n_0/2$. Thus $U_{n_0}(P,Q) = PU'_{n_0}(P,Q) = k\Box$ implies that $U'_{n_0}(P,Q) = k'\Box$, for $k'$ one of the finitely many divisors of $kn_0/2$. But $U_{n_0}(1, \frac{Q}{P^2}) = U'_{n_0}(1, \frac{Q}{P^2}) = P^{2-n_0}U'_{n_0}(P,Q)$, so that the curve $U_{n_0}(1,x) = k'y^2$ contains a rational point with $x = \frac{Q}{P^2}$. The polynomial $U_{n_0}(1,x)$, of degree $\lfloor\frac{n_0-1}{2}\rfloor$ has distinct roots: Schinzel [18], formula (70), page 58, gives this without proof, and for completeness we include a proof below. The genus of the hyperelliptic curve $U_{n_0}(1,x) = k''\Box$, $(k'' = k, k')$ is accordingly $\lfloor\frac{n_0-3}{4}\rfloor$, so at least equal to 2 when $n_0 \geq 11$, and the theorem will follow for $n_0 \geq 11$ from Faltings' proof of the Mordell Conjecture.

[1]The assumption $P > 0$ made in [15] does not affect the result, as $U_n(-P,Q)$ and $U_n(P,Q)$ differ at most in sign.

Write $u_{n_0}(x) = U_{n_0}(1,x)$, so that $u_{n_0}(x) = u_{n_0-1}(x) - xu_{n_0-2}(x)$, $u_0(x) = 0$, $u_1(x) = 1$. Solving the recurrence,

$$u_{n_0}(x) = \frac{1}{\sqrt{1-4x}}(z^{n_0} - \bar{z}^{n_0}), \quad z = \frac{1+\sqrt{1-4x}}{2}, \; \bar{z} = \frac{1-\sqrt{1-4x}}{2}.$$

It follows that the $\lfloor \frac{n_0-1}{2} \rfloor$ zeros of $u_{n_0}(x)$ occur where $z^{n_0} = \bar{z}^{n_0}$, i.e. where $z = \zeta\bar{z}$, $\zeta$ an $n_0$-th root of unity, i.e. where $\sqrt{1-4x} = \frac{\zeta-1}{\zeta+1}$. Thus the roots of $u_{n_0}(x)$ are given specifically by $x = \frac{\zeta}{(\zeta+1)^2}$, as $\zeta$ runs through the $n_0$-th roots of unity. Put $\zeta = e^{\frac{2\pi i j}{n_0}}$, $1 \le j \le \lfloor \frac{n_0-1}{2} \rfloor$; then the roots $x_j$ satisfy $x_j^{-1} = \zeta + 2 + \zeta^{-1} = 2 + 2\cos\frac{2\pi j}{n_0} = 4\cos^2\frac{\pi j}{n_0}$. Thus for $1 \le j \le \lfloor \frac{n_0-1}{2} \rfloor$ we obtain $\lfloor \frac{n_0-1}{2} \rfloor$ distinct roots, as required.

It remains to deal with $8 \le n_0 \le 10$, when the corresponding equations define a curve of genus 1 in weighted projective space.

Case $n_0 = 8$:

Now $U_8 = k\square$ implies

$$P = \lambda_1 a^2, \quad P^2 - 2Q = \lambda_2 b^2, \quad P^4 - 4P^2Q + 2Q^2 = \lambda_3 c^2, \tag{5}$$

with finitely many possibilities for squarefree $\lambda_i \in \mathbb{Z}$, $i = 1, 2, 3$. Eliminating $P$ and $Q$ at (5),

$$-\lambda_1^4 a^8 + 2\lambda_1^2\lambda_2 a^4 b^2 + \lambda_2^2 b^4 = 2\lambda_3 c^2, \tag{6}$$

and there may of course be infinitely many solutions to this equation: if the elliptic curve

$$-\lambda_1^4 A^4 + 2\lambda_1^2\lambda_2 A^2 B^2 + \lambda_2^2 B^4 = 2\lambda_3 C^2 \tag{7}$$

has positive rank, then a point $(A, B, C)$ leads to the point on (6) given by $(a, b, c) = (A, AB, A^2C)$. However, with the demand that $(a, b) = 1$, we show that there can only be finitely many solutions, as follows.

Let $(U_0, V_0, W_0)$ be an integer point on the underlying quadric

$$-U^2 + 2UV + V^2 = 2\lambda_3 W^2; \tag{8}$$

then by standard arguments (see, for example, formulas (20) in section 61 of [11]), the integer primitive solutions of the quadric (8) are of type

$$\begin{aligned}(U, V, W) \;=\; & (((U_0 - 2V_0)u^2 - 2V_0 uv + U_0 v^2)/\Delta, \\ & (-V_0 u^2 + 2U_0 uv - (2U_0 + V_0)v^2)/\Delta, \\ & W_0(u^2 - 2uv - v^2)/\Delta),\end{aligned}$$

where $(u,v)=1$ and $\Delta$ lies in a finite set. We thus have finitely many possibilities

$$\lambda_1^2 a^4 = ((U_0-2V_0)u^2 - 2V_0uv + U_0v^2)/\Delta, \qquad (9)$$

$$\lambda_2 b^2 = (-V_0u^2 + 2U_0uv - (2U_0+V_0)v^2)/\Delta, \qquad (10)$$

for some such $(u,v)$ and $\Delta$.

The $u,v$-quadratic at (10) is non-singular (its discriminant is $4(U_0^2-2U_0V_0-V_0^2)/\Delta^2 = -8\lambda_3W_0^2/\Delta^2 \neq 0$), and assuming it is rationally solvable, we can similarly parametrize the points giving finitely many possibilities

$$(u,v,b) = (q_1(r,s), q_2(r,s), q_3(r,s))$$

for quadratics $q_i$ and coprime $r,s$. Note that we cannot have $q_1(1,0)=q_2(1,0)=0$, otherwise $q_3(1,0)=0$ (from (10)) and the determinant of the coefficients of the quadratic forms $q_1$,$q_2$,$q_3$ is zero, impossible. Substituting this parametrization into (9), we obtain finitely many quartic curves of type

$$c\,a^4 = \text{quartic in } r,s,$$

namely

$$\Delta\lambda_1^2 a^4 = q(r,s) = (U_0-2V_0)q_1(r,s)^2 - 2V_0q_1(r,s)q_2(r,s) + U_0q_2(r,s)^2. \qquad (11)$$

We claim these quartic curves are non-singular, and hence of genus 3. First, the quartic $q(r,s)$ is not identically zero. For $U_0-2V_0$ and $U_0$ cannot both be zero, so suppose without loss of generality that $U_0-2V_0 \neq 0$. If $q(r,s)$ is identically zero, then $q_1(\lambda,1)/q_2(\lambda,1)$ is a root of the trinomial $(U_0-2V_0)X^2-2V_0X+U_0$ for more than two distinct values of $\lambda$, so that $q_1(r,s)=k_1q_2(r,s)$ for some constant $k_1$. Substitution into (10) implies then $q_3(r,s)=k_3q_2(r,s)$ for some constant $k_3$. This contradicts the fact that the determinant of the coefficients of the quadratic forms $q_1,q_2,q_3$ is non-zero.

Now singularity of the curve (11) occurs precisely when the quartic $q(r,s)$ contains a repeated factor, so we have two cases to consider.

Case I: $q(r,s)$ contains a repeated linear factor. By linear transformation on $r,s$ we can suppose the quartic is of type

$$q(r,s) = s^2 \text{ quadratic}(r,s). \qquad (12)$$

The coefficient of $r^4$ in $q(r,s)$ at (11) is $(U_0-2V_0)q_1(1,0)^2 - 2V_0q_1(1,0)q_2(1,0) + U_0q_2(1,0)^2$ which must therefore be zero. Since $q_1(1,0)$ and $q_2(1,0)$ are not both zero, then necessarily the discriminant $4V_0^2 - 4(U_0 - 2V_0)U_0 = 8\lambda_3 W_0^2$ must be a perfect square, forcing $\lambda_3 = 2$. Thus $P$ is even, $Q$ is odd, and $\lambda_2$ is also even. The quadric (8) has become

$$-U^2 + 2UV + V^2 = 4W^2,$$

so that with $(U_0, V_0, W_0) = (0, -2, -1)$ the equations at (9) and (10 ) are now:

$$\begin{aligned} \lambda_1^2 \Delta a^4 &= 4u(u+v), \\ \frac{\lambda_2}{2}\Delta b^2 &= u^2 + v^2. \end{aligned}$$

Rational solvability of the second quadric implies $\frac{\lambda_2}{2}\Delta = m^2 + n^2$, say, so that the parametrization $(u, v, b) = (q_1(r,s), q_2(r,s), q_3(r,s))$ is given by

$$\begin{aligned} q_1(r,s) &= (mr^2 + 2nrs - ms^2)/\delta, \\ q_2(r,s) &= (nr^2 - 2mrs - ns^2)/\delta, \\ q_3(r,s) &= (r^2 + s^2)/\delta \end{aligned}$$

for finitely many choices of $\delta$. The quartic at (11) becomes:

$$\lambda_1^2 \delta^2 \Delta a^4 = 4(mr^2 + 2nrs - ms^2)((m+n)r^2 - 2(m-n)rs - (m+n)s^2),$$

and the discriminant of the right-hand side is $512(m^2 + n^2)^6 \neq 0$, so the quartic cannot have a repeated root.

Case II: $q(r,s)$ at (11) contains a squared quadratic factor. But now the curve (11) is of genus 0, allowing $a, r, s$, hence $u, v, b$, to be written as polynomials in a single variable, so that the curve (7) is rationally parametrizable and of genus 0, which it is not. Thus the quartic at (11) is non-singular, and (11) defines a curve of genus 3. The finitely many such curves mean the original curve at (6) can have only finitely many points with coprime $(a, b)$.

Case $n_0 = 9$:

When $n_0 = 9$, then $U_9 = k\square$ implies $P^2 - Q = \lambda_1 R^2$, $P^6 - 6P^4Q + 9P^2Q^2 - Q^3 = \lambda_2\square$, with only finitely many possibilities for the $\lambda_i \in \mathbb{Z}$, $i = 1, 2$. Eliminating $Q$,

$$3P^6 - 9\lambda_1 P^4 R^2 + 6\lambda_1^2 P^2 R^4 + \lambda_1^3 R^6 = \lambda_2 S^2,$$

and the sextic has discriminant $-2^6 3^9 \lambda_1^{15}$, so is non-singular. Thus $x = P/R$ gives a point on a curve hyperelliptic of genus 2, and we are done.

Case $n_0 = 10$:

We have

$$P = \lambda_1 a^2, \quad P^4 - 3P^2Q + Q^2 = \lambda_2 b^2, \quad P^4 - 5P^2Q + 5Q^2 = \lambda_3 c^2, \qquad (13)$$

for finitely many squarefree $\lambda_i \in \mathbb{Z}$, $i = 1, 2, 3$. Thus $(x, y) = (P^2, Q)$ satisfies

$$x(x^2 - 3xy + y^2)(x^2 - 5xy + 5y^2) = \lambda_2 \lambda_3 \Box,$$

and finiteness of the number of solutions follows from Theorem 2 of Darmon and Granville [7].

**Remark**. Theorem 2.1 is not true for $n_0 = 2, \ldots, 7$; in fact, it is easy to check that for each fixed $n_0$ in this range, $U_{n_0}(P, Q) = \Box$ has infinitely many solutions $(P, Q)$ with $\gcd(P, Q) = 1$. See Section 2 of [3].

Assume now that we have computed explicitly all pairs $(P, Q)$ satisfying (2) and (4). Then, for any such specific pair $(P, Q) = (P_0, Q_0)$ we must find all positive indices $n = n_0 m$ with $\mathcal{P}(m) \subseteq \mathcal{S}$, satisfying

$$U_n(P_0, Q_0) = k\Box\,. \qquad (14)$$

In section 3 we show that, for a fixed pair $(P_0, Q_0)$ as above, it is possible, at least in principle, to solve explicitly equation (14) in the unknown $n$ without any restriction on the prime divisors of $n$.

## 3 Terms of Lucas sequences being almost squares

In this section we will show that for fixed non-zero integer $k$ and parameters $P, Q$ which satisfy (2), the problem of finding all $n$ for which

$$U_n(P, Q) = k\Box\,, \qquad (15)$$

leads to finitely many quartic Thue-Mahler equations, where the subscript $n$ will appear in the exponents of the prime numbers in the right-hand side of these equations. Thue-Mahler equations, in general, can be explicitly solved, at least in principle, by the method explained in great detail in [21]. It is of some interest to note at this point that for our special purposes we do not need the complete solution of the Thue-Mahler equations but only a small upper bound for the unknown exponents of the right-hand side; we make a short discussion of this issue at the beginning of section 6.

By section 1 of [15] we have

$$V_n(P,Q)^2 - DU_n(P,Q)^2 = 4Q^n \tag{16}$$

and $\gcd(U_n(P,Q), V_n(P,Q)) = 1$ or 2.

Suppose first that $n$ is even, $n = 2m$, say. Then $(2Q^m, U_{2m}(P,Q), V_{2m}(P,Q))$ is an integer solution of the equation

$$X^2 + DY^2 - Z^2 = 0\,, \quad (X,Y,Z) = 1 \text{ or } 2\,. \tag{17}$$

A special primitive solution of (17) is $(X,Y,Z) = (-1,0,1)$, and therefore all primitive solutions are given by

$$(X,Y,Z) = (\frac{S^2 - DT^2}{\Delta}, \frac{2ST}{\Delta}, \pm\frac{S^2 + DT^2}{\Delta})\,, \tag{18}$$

where $\gcd(S,T) = 1$ and $\Delta$ is, up to sign, the gcd of the numerators of the fractions in the right-hand side. It follows that $\Delta$ must divide the determinant of the coefficients of the three quadratic forms in $S, T$ appearing in these numerators: therefore $\Delta|4D$. Now $2^{\nu+1}ST = \Delta U_n = k\Delta\Box$ and $2^\nu(S^2 - DT^2) = 2\Delta Q^m$, where $\nu = 0$ or 1, according as $\gcd(U_{2m}(P,Q), V_{2m}(P,Q)) = 1$ or 2. Since $\gcd(S,T) = 1$, it follows from the first equation that $S = k_1u^2$, $T = k_2v^2$, where $\gcd(u,v) = 1$ and $k_1, k_2$ are integers that can be easily computed explicitly by means of $k$ and $\Delta$. Substitution of these expressions for $S$ and $T$ into the second equation gives

$$k_1^2u^4 - Dk_2^2v^4 = \ \Delta Q^m \ \text{ or } \ 2\Delta Q^m\,.$$

Obviously, the polynomial $k_1t^4 - Dk_2^2$ has distinct complex roots, and therefore the above equation is a Thue-Mahler equation.

Suppose second that $n$ is odd, $n = 2m+1$, say. Then $(2Q^m, U_{2m+1}(P,Q), V_{2m+1}(P,Q))$ is a solution $(X,Y,Z)$ of

$$QX^2 + DY^2 - Z^2 = 0\,, \quad (X,Y,Z) = 1 \text{ or } 2\,. \tag{19}$$

A special primitive solution of (19) is $(X,Y,Z) = (-2,-1,P)$, and therefore all primitive solutions are given by

$$(X,Y,Z) = (\frac{2QS^2 + 2DST - 2DT^2}{\Delta}, \frac{-QS^2 + 4QST + DT^2}{\Delta}, \pm P\frac{QS^2 + DT^2}{\Delta})\,, \tag{20}$$

where, as before, $\Delta$ (up to sign) is the gcd of the numerators of the fractions in the right-hand side. Hence $\Delta | 2P^3QD$.

Now $2^\nu(-QS^2+4QST+DT^2) = \Delta U_n = k\Delta\Box$ and $2^\nu(QS^2+DST-DT^2) = \Delta Q^m$, where $\gcd(S,T) = 1$ or 2 and, as before, $\nu = 0$ or 1, according as $\gcd(U_{2m}(P,Q), V_{2m}(P,Q)) =$ 1 or 2. The first equation leads to

$$k\Delta X^2 - 2^\nu P^2 Y^2 + 2^\nu Q Z^2 = 0\,, \tag{21}$$

where $Y = T$ and $Z = S - 2T$, so that $\gcd(X,Y,Z) = 1$ or 2; and the second equation then becomes

$$2^\nu(QZ^2 + P^2ZY + P^2Y^2) = \Delta Q^m\,. \tag{22}$$

Equation (21) implies an equation $Ax^2+By^2+Cz^2 = 0$, where $A, B, C$ are non-zero integers and $ABC$ is square-free[2]. For the non-trivial solvability of this equation there is the classical Legendre criterion, which is easily applied. Finding an actual solution is a much more difficult and computationally interesting problem, especially if the size of $A, B, C$ is large; for the discussion of this problem we refer the reader to [6]. Here, we assume that we know an integer solution $(x_0, y_0, z_0)$ of (21) with $\gcd(x_0,y_0,z_0) = 1$ and $z_0 \neq 0$. Then all integer solutions to (21) with $\gcd(X,Y,Z) =$ 1 are given by the formulas

$$\begin{aligned} \delta X &= -k\Delta x_0 \cdot u^2 - 2^{\nu+1}P^2 y_0 \cdot uv + 2^\nu P^2 x_0 \cdot v^2 \\ \delta Y &= k\Delta y_0 \cdot u^2 - 2k\Delta x_0 \cdot uv - 2^\nu P^2 y_0 \cdot v^2 \\ \delta Z &= \pm z_0(k\Delta \cdot u^2 + 2^\nu P^2 \cdot v^2) \end{aligned}$$

where $\gcd(u,v) = 1$ and $\delta > 0$ is the gcd of the numbers in the three right-hand sides; $\delta$ must divide the determinant of the coefficient matrix of the three quadratic forms in the right-hand sides, which is equal to $2^{2\nu+2}k\Delta P^2 Q z_1^3$.

Substitution of the above expressions for $Y, Z$ into equation (22) gives a quartic form in $u, v$ being equal to $\Delta Q^m$. The discriminant of the corresponding quartic polynomial is equal to $-2^{8+6\nu}(DQ)^3(k\Delta)^6(Pz_0)^{12} \neq 0$, so that we have arrived at a quartic Thue-Mahler equation.

[2]In general, the implied equation is not equivalent, from the point of view of integer solutions, to (21).

# 4 An example: solution of $U_n(P,Q) = \pm\Box$ or $\pm 2\Box$ for $n$ a power of 2

As an application of the discussion in the previous sections we will solve the equations

$$U_{2^e}(P,Q) = \pm\Box, \pm 2\Box\ ,$$

in the unknowns $e, P, Q$, with $e \geq 3$ and $P, Q$ satisfying (2). In the notation of Problem 1, $n_0 = 8$, $m = 2^{e-3}$, $\mathcal{S} = \{2\}$ and $k = \pm 1$ or $\pm 2$. Relation (4) now becomes

$$U_8(P,Q) = \pm 2^{i_1}\Box\,, \quad i_1 \in \{0,1\}\ .$$

In this connection we have the following result.

**Theorem 4.1.** *Let $k_0|2$. The only solutions $(P,Q)$ of $U_8(P,Q) = k_0\Box$ satisfying (2) are $(P,Q) = (k_0,-4), (4k_0,-17)$ when $k_0 = \pm 1$, while no solutions exist when $k_0 = \pm 2$.*

*Proof.* See section 5. $\Box$

Theorem 4.1 implies $i_1 = 0$ and $(P,Q) = (\pm 1,-4), (\pm 4,-17)$. Since for $n$ even, $U_n(P,Q) = -U_n(-P,Q)$, without loss of generality it remains only to find explicitly all $n = 2^e$, with $e \geq 3$, such that

$$U_n(\pm 1,-4) = k\Box, \quad k = 1, 2 \tag{23}$$

$$U_n(\pm 4,-17) = k\Box, \quad k = 1, 2\ . \tag{24}$$

Consider the case $k = 1$ at (23) and write for simplicity $U_n, V_n$ instead of $U_n(\pm 1,-4)$, $V_n(\pm 1,-4)$. Relation (16) now reads $(2^{n+1})^2 + 17U_n^2 = V_n^2$, leading to an equation of the form $X^2 + 17Y^2 = Z^2$, where $X$ is even and $\gcd(X, 17Y) = 1$. This is an equation as at (17) and from the formulas (18) we obtain $\pm 2^{n+1} = (17T^2 - S^2)/2$ and $U_n = ST$, where $ST$ is odd and $\gcd(S, 17T) = 1$. Since $U_n = \Box$, we are finally led to

$$u^4 - 17v^4 = \pm 2^{n+2}\,, \quad uv \text{ odd}\,. \tag{25}$$

From Proposition 6.1, it follows that $n = 8$.

Second, consider the case $k = 2$ at (23). Since $U_n(\pm 1,-4)$ is always odd, this case is impossible.

Third, consider the case $k = 1$ at (24). As before, write $U_n, V_n$ instead of $U_n(\pm 4, -17)$, $V_n(\pm 4, -17)$. Relation (16) becomes $(17^{n/2})^2 + 21U_n^2 = (V_n/2)^2$. It is easy to check that, except for $U_0$, every $U_n$ and every $V_n$ is indivisible by 17. Moreover, for $n$ even, $U_n$ is even. Therefore we are led to an equation (17) with $D = 21$, in which $X$ is odd, $Y$ is even and $\gcd(17X, Y) = 1$. By the formulas (18) or otherwise, and in view of the fact that $U_n = \square$, we obtain the equations

$$\pm 17^{n/2} = 21b^4 - 4a^4\,, \quad \pm 17^{n/2} = 84a^4 - b^4\,,$$

$$\pm 17^{n/2} = 7b^4 - 12a^4\,, \quad \pm 17^{n/2} = 28a^4 - 3b^4\,,$$

in which $b$ is odd, $\gcd(a, b) = 1$ and $ab \not\equiv 0 \bmod 17$. The last two equations are impossible mod 17; and in the first equation we see mod 4 that the plus sign holds, whereas mod 3 the minus sign holds (since $n/2$ is even), impossible. In the second equation the minus sign must hold, and we finally get, on putting $(u, v) = (b, a)$,

$$u^4 - 84v^4 = 17^{n/2} \tag{26}$$

where $uv \not\equiv 0 \bmod 17$. From Proposition 7.1 it follows that $n = 8$.

Last, consider the case $k = 2$ at (24). Arguing as as above, we are led to the equations $\pm 17^{n/2} = 7u^4 - 3v^4$, $\pm 17^{n/2} = u^4 - 21v^4$, of which the former is impossible mod 17. In the latter, a congruence mod 3 shows that the plus sign must hold, resulting in the equation

$$u^4 - 21v^4 = 17^{n/2}. \tag{27}$$

This equation is very similar to (26). Indeed, although the quartic fields corresponding to $x^4 - 84$ and $x^4 - 21$ are not isomorphic, both have class number two, and in both fields 17 splits into four prime ideals, of which two are principal, and two are non-principal. Solving the equation $u^4 - 21v^4 = 17^z$ using the same techniques as illustrated in Section 7 shows that the only solution with $n > 0$ of (27) is $(u, v, n) = (5, 2, 4)$.

We have now proved the following result.

**Theorem 4.2.** *The only solutions $(2^e, P, Q)$ with $e \geq 3$ of the equation*

$$U_{2^e}(P, Q) = \pm\square$$

*under the constraint (2) are* $(8, \pm 1, -4)$ *and* $(8, \pm 4, -17)$*. There are no solutions with* $e \geq 3$ *to the equation*

$$U_{2^e}(P, Q) = \pm 2\square.$$

## 5 Proof of Theorem 4.1

The demand that $U_8(P, Q) = k_0\square$, $k_0|2$, is that $P(P^2-2Q)(P^4-4P^2Q+2Q^2) = k_0\square$, and $P$, $Q$ must satisfy the constraint (2). Necessarily there exist integers $a, b, c$, coprime in pairs, satisfying $P = \delta_1 a^2$, $P^2 - 2Q = \delta_2 b^2$, $P^4 - 4P^2Q + 2Q^2 = \delta_3 c^2$, with $\delta_i|2$, $i = 1, 2, 3$, and $k_0 \equiv \delta_1\delta_2\delta_3 \bmod \mathbb{Q}^{*2}$. Eliminating $P$, $Q$,

$$-\delta_1^4 a^8 + 2\delta_1^2\delta_2 a^4 b^2 + \delta_2^2 b^4 = 2\delta_3 c^2. \tag{28}$$

It is straightforward to check that (28) is everywhere locally solvable only in the following instances:

$$(\delta_1, \delta_2, \delta_3) = (\pm 1, \pm 2, 2), (\pm 1, 1, 1), (\pm 1, -1, -1), (\pm 2, \pm 2, 2). \tag{29}$$

With the first triple at (29), $a, b, c$ satisfy $-a^8 \pm 4a^4b^2 + 4b^4 = 4c^2$, so that $a$ is even, and $(A, B, C) = (a/2, b, c)$ satisfy $B$, $C$ odd, and

$$-64A^8 \pm 16A^4B^2 + B^4 = C^2; \tag{30}$$

for the second triple at (29),

$$-a^8 + 2a^4b^2 + b^4 = 2c^2; \tag{31}$$

for the third triple at (29),

$$-a^8 - 2a^4b^2 + b^4 = -2c^2; \tag{32}$$

and for the fourth triple at (29),

$$-4a^8 \pm 4a^4b^2 + b^4 = c^2. \tag{33}$$

It suffices to find all integer solutions to the equations (30)-(33), with coprime $(A, B)$ and $(a, b)$, as appropriate. The equations (30), (31), (32) have been treated in Bremner & Tzanakis [3]. The only solutions to (30) with $+$ sign, i.e. for $\delta_2 = 2$, having $(A, B) = 1$ and $B$ odd, are $(\pm A, \pm B, \pm C) = (0, 1, 1), (1, 5, 31)$; and with $-$ sign, i.e. for $\delta_2 = -2$, the only solutions are $(\pm A, \pm B, \pm C) = (0, 1, 1)$. The sole resulting $(P, Q)$ satisfying (2) arises from $(1, 5, 31)$, with $(P, Q) = (4\delta_1, -17) = (4k_0, -17)$, as required. The only solutions to (31), (32), are respectively $(\pm a, \pm b, \pm c) = (1, 1, 1), (1, 3, 7)$; and $(1, 1, 1)$. The sole resulting $(P, Q)$ satisfying (2) arises from $(1, 3, 7)$, with $(P, Q) = (\delta_1, -4) = (k_0, -4)$. Consider finally (33), where a solution $(a, b, c)$ implies a solution $(A, B, C) = (a, 2b, 4c)$ of (30). Thus (33) may be treated similarly to the equation (30), where we instead seek solutions with $B$ *even.* The modifications necessary are very minor, and the only solutions to (33) with $+$ sign are $(\pm a, \pm b, \pm c) = (0, 1, 1), (1, 1, 1)$; and the only solutions with $-$ sign are $(\pm a, \pm b, \pm c) = (0, 1, 1)$. None of these leads to suitable $P, Q$.

# 6 Solution of $u^4 - 17v^4 = \pm 2^{z+2}$

In this section we prove the following:

**Proposition 6.1.** *The only solutions of the Diophantine equation* $u^4 - 17v^4 = \pm 2^{z+2}$ *with* $uv$ *odd and* $z \geq 1$ *are given by* $1^4 - 17 \cdot 1^4 = -2^4$, $3^4 - 17 \cdot 1^4 = 2^6$ *and* $7^4 - 17 \cdot 3^4 = 2^{10}$.

Note that, for the needs of the proof of Theorem 4.2 a considerably slighter result suffices, namely, that $z$ is bounded by a small upper bound, say, of the size of 1000. For, in this case, we can very easily check all $n$ that are powers of 2, for which $U_n(\pm 1, -4)$ has the required shape. We decided, however, to solve completely the Thue-Mahler equation of Proposition 6.1, as this can be accomplished, at least in this special example, with very little extra cost. Note that the curves $x^4 - 17y^4 = cz^4$, $\pm c = 1, 2, 4, 8$, of genus 3, all have Jacobian of rank at least 5, and direct Chabauty arguments are inapplicable.

## 6.1 Preliminaries to the solution

We work in the field $K = \mathbb{Q}(\theta)$, where $\theta^4 = 17$. Using PARI [13] we obtained the following information: The class-number of $K$ is 2, an integral basis is $1, \theta, \psi = (1 + \theta^2)/2, \omega = (1 + \theta + \theta^2 + \theta^3)/4$ and a pair of fundamental units is $\epsilon_1 = 2 - \theta, \epsilon_2 = 4 - \theta^2$, with respective norms -1 and +1. The ideal factorization of 2 is

$$\langle 2 \rangle = \mathfrak{p}_1 \mathfrak{p}_2 \mathfrak{p}_3^2 \,, \; \mathfrak{p}_1 = \langle 2, -1 + \omega \rangle \,, \; \mathfrak{p}_2 = \langle 2, -1 + \psi + \omega \rangle \,, \; \mathfrak{p}_3 = \langle 2, 1 - \theta + \psi \rangle$$

and the ideal-class of $\mathfrak{p}_2$ generates the ideal-class group of $K$.

Factorization of (25) in $K$ gives

$$(u - v\theta)(u + v\theta)(u^2 + v^2\theta^2) = \pm 2^z \,. \tag{34}$$

We have the ideal factorizations

$$\langle 1-\theta\rangle = \mathfrak{p}_1^2\mathfrak{p}_2\mathfrak{p}_3\,,\quad \langle 1+\theta\rangle = \mathfrak{p}_1\mathfrak{p}_2^2\mathfrak{p}_3\,,\quad \langle 1-\theta^2\rangle = \mathfrak{p}_1^3\mathfrak{p}_2^3\mathfrak{p}_3^2\,,\qquad \langle 1+\theta^2\rangle = \mathfrak{p}_1\mathfrak{p}_2\mathfrak{p}_3^6\,.$$

Without loss of generality we may assume that $u \equiv v \bmod 4$. Then, with the aid of the above relations it is easy to see that the ideal $\langle u - v\theta\rangle$ is divisible by $\mathfrak{p}_i$ for $i = 1, 2, 3$, but is not divisible by a higher power of $\mathfrak{p}_i$ for $i = 2, 3$. Also, each ideal $\langle u + v\theta\rangle$ and $\langle u^2 + v^2\theta^2\rangle$ is divisible by $\mathfrak{p}_1$ but by no higher power of $\mathfrak{p}_1$. This, in combination with the ideal equation obtained from (34) leads to the ideal equation[3]

$$\langle u - v\theta\rangle = \mathfrak{p}_1^{z-2}\mathfrak{p}_2\mathfrak{p}_3\,. \tag{35}$$

The ideals $\mathfrak{p}_2$ and $\mathfrak{p}_3$ belong to the same ideal-class and $\mathfrak{p}_2^2$ is principal. Therefore, the above ideal equation implies that $z$ must be even. We put

$$z - 2 = 2n_1\,.$$

With the aid of PARI we see that

$$\mathfrak{p}_2\mathfrak{p}_3 = \langle 1 + \theta + \omega\rangle\,,\quad \mathfrak{p}_1^2 = \langle -2 + \theta - 2\psi + \omega\rangle\,,$$

so that (35) finally leads to

$$u - v\theta = \pm\alpha\epsilon_1^{a_1}\epsilon_2^{a_2}\pi_1^{n_1}\,, \tag{36}$$

where

$$\alpha = 1+\theta+\omega = \frac{1}{4}(5+5\theta+\theta^2+\theta^3)\,,\quad \pi_1 = -2+\theta-2\psi+\omega = \frac{1}{4}(-11+5\theta-3\theta^2+\theta^3).$$

**Notations and conventions in what follows**. The conjugates of the typical algebraic number $\gamma$, say, considered either as $p$-adic numbers ($p = 2$), belonging to a finite extension of $\mathbb{Q}_2$, or as complex numbers, will be denoted, in both cases, by $\gamma^{(i)}$, $i = 1, 2, \ldots$. However, there is no fear of confusion, as it will be absolutely clear from the text whether we work in $\mathbb{C}$ or in an extension of $\mathbb{Q}_2$.

The 2-adic integer $b_0 + 2b_1 + 2^2b_2 + \cdots$, where the $b_i$'s are binary digits, is written in the form $0.b_0b_1b_2\ldots$. Also, $\mathrm{ord}_2$ will denote the 2-adic additive valuation, which is defined on any finite extension of $\mathbb{Q}_2$ and extends the usual 2-adic additive valuation of $\mathbb{Q}$. For a quick practical survey we refer to section 4 of [21].

[3] On considering equation (25) mod 8 we see that $z \geq 3$.

As usual, $i$ stands for the root of the polynomial $t^2+1 \in \mathbb{Q}[t]$ and $\imath$ will denote its "value", complex or 2-adic, as the case may be.

Similarly, $\vartheta$ will denote the value of $\theta$ which, in the "real context" will mean the number $\sqrt{17} = 4.123105\ldots \in \mathbb{R}$ and in the "2-adic context" the number $\sqrt{17} = 0.101101011011\ldots \in \mathbb{Q}_2$.

## 6.2 A large upper bound for $\max\{n_1, |a_1|, |a_2|\}$

Following the notation of [21] we put $A = \max\{|a_1|, |a_2|\}$ and $H = \max\{A, n_1\}$. The various $c$ constants with subscripts, which we will occasionally mention below always agree with those in [21].

We set

$$\theta^{(1)} = \vartheta\,, \quad \theta^{(2)} = -\vartheta\,, \quad \theta^{(3)} = \imath\vartheta\,, \quad \theta^{(4)} = -\imath\vartheta\,. \tag{37}$$

and

$$\lambda = \delta\left(\frac{\pi_1^{(4)}}{\pi_1^{(3)}}\right)^{n_1}\left(\frac{\epsilon_1^{(4)}}{\epsilon_1^{(3)}}\right)^{a_1}\left(\frac{\epsilon_2^{(4)}}{\epsilon_2^{(3)}}\right)^{a_2} - 1 = \delta\left(\frac{\pi_1^{(4)}}{\pi_1^{(3)}}\right)^{n_1}\left(\frac{\epsilon_1^{(4)}}{\epsilon_1^{(3)}}\right)^{a_1} - 1$$

(note that $\epsilon_2^{(3)} = \epsilon_2^{(4)}$), where $\delta = \dfrac{\theta^{(1)}-\theta^{(3)}}{\theta^{(1)}-\theta^{(4)}}\cdot\dfrac{\alpha^{(4)}}{\alpha^{(3)}}$. Below we will view $\lambda$ both as a 2-adic and a complex number. In the first case we consider $\Lambda_1 = \log_2(\lambda+1)$, where $\log_2$ stands for the logarithm in the $p$-adic sense (see section 12 of [21]) with $p=2$; and in the second case we consider $\Lambda_0 = \imath^{-1}\mathrm{Log}\,(1+\lambda)$, where Log, here and everywhere in this paper, stands for the principal branch of the complex logarithmic function.

It is worth noticing at this point that $\epsilon_2$, hence $a_2$ as well, "disappeared" from $\lambda$ and consequently from $\Lambda_0$ and $\Lambda_1$. However, "hidden" in $A$, $a_2$ indirectly appears in the absolutely crucial relation

$$0 < |\Lambda_0| < 1.02 c_{21} e^{-c_{16}A} \tag{38}$$

(see relation (27) in [21]). Moreover, it should be stressed that both $\epsilon_1$ and $\epsilon_2$ play a significant role in the computations of the numerous positive constants[4] leading to the upper bound for $H$.

[4]Such a constant is, in the notation of [21], $c_{15}$ and all constants depending on it.

The numbers (complex or 2-adic) appearing in $\lambda$ are

$$\begin{aligned}
\delta &= \frac{1}{8}(\vartheta - \vartheta^3 + (9-\vartheta^2)\imath) \\
\frac{\pi_1^{(4)}}{\pi_1^{(3)}} &= \frac{1}{16}(33 - 9\vartheta^2 + (13\vartheta - 5\vartheta^3)\imath) = -\delta^{-2} \\
\frac{\epsilon_1^{(4)}}{\epsilon_1^{(3)}} &= -33 + 8\vartheta^2 - (16\vartheta - 4\vartheta^3)\imath := \beta \,,
\end{aligned}$$

hence,

$$\lambda = (-1)^{n_1}\delta^{1-2n_1}\beta^{a_1} - 1 \,. \tag{39}$$

First we work $p$-adically, with $p = 2$. Following [21][5] we see that, for $n_1$ positive, $\mathrm{ord}_2(\lambda) = 2n_1$. On the other hand, with the aid of the *Theory of Linear Forms in p-adic Logarithms* applied to $\lambda$, which we write as $(\pm\delta)^{1-2n_1}\beta^{a_1} - 1$, we can bound from above $\mathrm{ord}_2(\lambda)$ in terms of the logarithm of the maximum absolute value of the integer unknowns appearing in (39), i. e. in terms of $\log\max(|a_1|, 2n_1 - 1)$, hence in terms of $\log H$. Actually, using a very useful theorem of Kunrui Yu [24][6] we compute constants $c_{13}, c_{14}$, such that

$$n_1 \le c_{13}(\log H + c_{14}) \,, \tag{40}$$

where the large constant is always $c_{13}$. In our example, this is of the size of $10^{19}$.

Next we work in $\mathbb{C}$. Now each of the complex numbers $\beta$ and $\delta$ in (39) is the ratio of two complex-conjugate numbers, and therefore the Log is at most $\pi$ in absolute value. On the other hand, in general, $\mathrm{Log}(z_1 + z_2) \equiv \mathrm{Log}\, z_1 + \mathrm{Log}\, z_2 \pmod{2\pi}$, therefore

$$\begin{aligned}
\imath\Lambda_0 =& \mathrm{Log}\,\delta + n_1\mathrm{Log}(-1) - 2n_1\mathrm{Log}\,\delta + a_1\mathrm{Log}\,\beta + 2a_0\mathrm{Log}(-1) \\
=& (2a_0 + n_1)\mathrm{Log}(-1) + a_1\mathrm{Log}\,\beta + (1-2n_1)\mathrm{Log}\,\delta \,,
\end{aligned} \tag{41}$$

so that a new unknown integer $a_0$ makes its appearance, satisfying $|a_0| \le 2H + 1$.

With the aid of the *Theory of Linear Forms in Real/Complex Logarithms* we can bound from below $\imath\Lambda_0$ in terms of $\log\max(|2a_0 + n_1|, |a_1|, 2n_1 - 1)$, hence in terms of $\log H$. We applied the corollary to the main theorem[7] of A. Baker and G. Wüstholz [1] to obtain $|\Lambda_0| > \exp\{-c_7(\log H + c_8')\}$, with explicit constants $c_7, c_8$

[5] cf. relation (13) therein.

[6] See also Appendix A2 in [21] for a version of Yu's theorem specially adapted to our needs.

[7] in fact, the application of this corollary as proposed in [22].

of which $c_7$ is the really large one; in our case this is of the size of $10^{27}$. Combined with (38) this gives a relation $A < C_3 \log H + C_4$, with explicit constants, and from (40) we have a relation $n_1 < C_1 \log H + C_2$. In view of the definition of $H$, these two inequalities combine to bound $H$ by a linear function of $\log H$, implying thus an upper bound for $H$. In our case this bound is of the size of $10^{29}$.

The rather long sequence of boring computations culminating in the upper bound for $H$, including the application of the theorem of Yu and that of Baker and Wüstholz, were performed "almost" automatically with the aid of MAPLE. We say "almost" because we computed the minimal polynomials of various algebraic numbers first, which in the sequel we inserted as (part of the) input of our MAPLE procedure. Our code is available upon request.

### 6.3 All solutions of equation (25)

Having obtained such large upper bounds for $n_1$ and $H$ we must reduce them in order to be able to solve our equation *explicitly.* We follow the reduction strategy as explained in great detail in sections 13 through 16 of [21]. The main tool for the reduction is the LLL algorithm, which is applied to appropriate two or three dimensional lattices generated by basis vectors with rational integer coordinates. In one instance, among the basis coordinates there appears an integer which is a 2-adic approximation of $\log_2 \beta / \log_2 \delta$ to a precision $2^m$, with $m$ sufficiently large, depending on the size of the upper bound for $n_1$. In our case it has not been necessary to take $m$ larger than 170. This is the *p-adic reduction step.* In the other instance, among the basis coordinates there appear the integer parts of the real numbers $C \mathrm{Arg}\, \delta, C \mathrm{Arg}\, \beta$ and $C\pi$, where $C$ is a sufficiently large integer, depending on the size of the upper bound for $H$. In our case, the larger value for $C$ that we needed was $10^{65}$. This is the *real reduction step.*

One starts with the initial large upper bounds for $n_1$ and $H$, say $N_0$ and $H_0$, respectively. The $p$-adic reduction step is applied first and reduces $N_0$ to a considerably lower bound $N_1$ which, in our case, turned out to be 85. According to the remark we made at the beginning of section 6, if we are interested just in proving Theorem 4.2 and not in the Thue-Mahler equation itself, we can stop here and easily check for which $n = 8, 16, 32, 64$ it is true that $U_n(\pm 1, -4)$ is, up to sign, a square or two times a square. It turns out that only $n = 8$ is acceptable.

However, at a little more cost, we can solve completely our Thue-Mahler equation, as follows. After the above reduction of the upper bound of $n_1$, we apply the real reduction step taking into account the bounds $N_1$ and $H_0$. As a result, a reduction of $H_0$ to a considerably smaller bound $H_1$ (less than 600) is achieved. The $p$-adic reduction step is again applied, taking into account the bounds $N_1$ and $H_1$ and returning a smaller upper bound $N_2$ etc. The iterative application of the $p$-adic and real reduction goes on until no further essential improvement is achieved. In our case, three couples of iterations had as a result the bounds $n_1 \leq 10$ and $A \leq 32$.

Now we can run through all $(n_1, a_1, a_2)$ in the range $1 \leq n_1 \leq 10$, $-32 \leq a_i \leq 32$ to check for which of them it is true that the coefficients of both $\theta^2$ and $\theta^3$ in $\alpha\epsilon_1^{a_1}\epsilon_2^{a_2}\pi_1^{n_1}$ are zero; cf. (36). Although there is no problem in checking 42250 triads $(n_1, a_1, a_2)$, we can do something better and diminish the number of checks by a factor of 65, as follows. We observe that if $(n_1, a_1, a_2)$ is such that (36) holds, then

$$\begin{aligned}\frac{u - v\theta}{u + v\theta} &= \frac{\alpha^{(1)}}{\alpha^{(2)}}\left(\frac{\epsilon_1^{(1)}}{\epsilon_1^{(2)}}\right)^{a_1}\left(\frac{\pi_1^{(1)}}{\pi_1^{(2)}}\right)^{n_1} \\ &= \frac{5 + 5\theta + \theta^2 + \theta^3}{5 - 5\theta + \theta^2 - \theta^3}\left(\frac{2-\theta}{2+\theta}\right)^{a_1}\left(\frac{-11 + 5\theta - 3\theta^2 + \theta^3}{-11 - 5\theta - 3\theta^2 - \theta^3}\right)^{n_1}\end{aligned}$$

If we denote the right-hand side by $R(n_1, a_1)$, then

$$\frac{1}{1 - R(n_1, a_1)} = \frac{1}{2} + \frac{u}{34v}\theta^3$$

and we must check for which pairs $(n_1, a_1)$ the left-hand side has rational part $1/2$ and both coefficients of $\theta$ and $\theta^2$ zero. It turns out that the only solutions $n_1$ with corresponding $u/v$ are $(n_1, u/v) = (1, 1), (2, -3), (4, 7/3)$, which completes the proof of Proposition 6.1.

# 7 Solution of $u^4 - 84v^4 = 17^z$

In this section we prove the following:

**Proposition 7.1.** *Under the constraint* $uv \not\equiv 0 \bmod 17$, *the only solution to the Diophantine equation* $u^4 - 84v^4 = 17^z$ *is* $31^4 - 84 \cdot 10^4 = 17^4$.

A comment analogous to that made immediately after the announcement of Proposition 6.1 is still valid: A small upper bound for the exponent $z$ appearing in the above equations suffices for the needs of the proof of Theorem 4.2. We work in the field $K = \mathbb{Q}(\theta)$, where $\theta^4 = 84$. The class-number of $K$ is 2, an integral basis is $1, \theta, \omega, \theta\omega$, where $\omega = (2 + \theta^2)/4$, and a pair of fundamental units is $\epsilon_1 = -55 + 18\theta - 6\theta^2 + 2\theta^3$, $\epsilon_2 = 2 + \omega$. The ideal factorization of 2 is $\langle 2 \rangle = \mathfrak{p}_2^2$ and the ideal class $[\mathfrak{p}_2]$ generates the ideal class group. Also, $\langle 17\rangle$ is the product of four ideals $\mathfrak{p}_{17,i}$, $i=1,\ldots,4$, two of which are principal and two non-principal:

$$\mathfrak{p}_{17,1}=\langle\theta-\omega\rangle\,,\ \mathfrak{p}_{17,2}=\langle 17,-2+\theta\rangle\,,\ \mathfrak{p}_{17,3}=\langle 17,2+\theta\rangle\,,\ \mathfrak{p}_{17,4}=\langle\theta+\omega\rangle\,,$$

where

$$\mathfrak{p}_{17,2}=(-1+\frac{\theta}{2})\mathfrak{p}_2 \quad\text{and}\quad \mathfrak{p}_{17,2}^2=\langle 1-2\theta+2\omega\rangle\,,$$

$$\mathfrak{p}_{17,3}=(1+\frac{\theta}{2})\mathfrak{p}_2 \quad\text{and}\quad \mathfrak{p}_{17,3}^2=\langle 1+2\theta+2\omega\rangle\,.$$

Consider now $u^4-84v^4=17^z$, with $uv\not\equiv 0 \bmod 17$. From the ideal factorization $\langle u-v\theta\rangle\langle u+v\theta\rangle\langle u^2+v^2\theta^2\rangle=\langle 17\rangle^z$, in which the ideals at left are pairwise relatively prime, we see that $\langle u-v\theta\rangle=\mathfrak{p}_{17,i}^z$ for some $i=1,\ldots,4$, with $i=2,3$, being possible only if z is even. It follows that $\langle u-v\theta\rangle$ must be equal to one of the following ideals:

$$\langle\theta-\omega\rangle^z\,,\ \langle 1-2\theta+2\omega\rangle^{z/2}\,,\ \langle\theta+\omega\rangle^z\,,\ \langle 1+2\theta+2\omega\rangle^{z/2}\,.$$

The third and fourth possibility are essentially the same as the first and second one, respectively, in view of the automorphism $\theta\mapsto-\theta$. Moreover, on choosing appropriately the signs of $u,v$, we obtain from the first two ideal equations the following element equations

$$u-v\theta \;=\; \epsilon_1^{a_1}\epsilon_2^{a_2}(\theta-\omega)^{n_1}\,,\ n_1=z \tag{42}$$

$$u-v\theta \;=\; \epsilon_1^{a_1}\epsilon_2^{a_2}(1-2\theta+2\omega)^{n_1}\,,\ n_1=z/2\,. \tag{43}$$

We work in complete analogy with section 6.1. Now $p=17$ and $\vartheta$, the value of $\theta$, in the "real context" will mean the number $\sqrt{84}=9.16515139\ldots\in\mathbb{R}$ and in the "17-adic context" the number $\sqrt{84}=0.8(13)(12)4(13)(13)\ldots\in\mathbb{Q}_{17}$.

The symbol $i$ still stands for the root of the polynomial $t^2+1\in\mathbb{Q}[t]$ and $\imath$ will denote its "value", complex or 17-adic, as the case may be.

We put

$$\pi_1=\begin{cases}\theta-\omega & \text{in the case of equation (42)}\\ 1-2\theta+2\omega & \text{in the case of equation (43)}\end{cases}$$

and, in both cases, we keep the notations in (37) with $\theta,\epsilon_1,\epsilon_2,\pi_1$ having, of course, their new values. In analogy with section 6.1 we set

$$\lambda=\frac{\vartheta^{(1)}-\vartheta^{(3)}}{\vartheta^{(1)}-\vartheta^{(4)}}\left(\frac{\pi_1^{(4)}}{\pi_1^{(3)}}\right)^{n_1}\left(\frac{\epsilon_1^{(4)}}{\epsilon_1^{(3)}}\right)^{a_1}\left(\frac{\epsilon_2^{(4)}}{\epsilon_2^{(3)}}\right)^{a_2}-1=-\imath\left(\frac{\pi_1^{(4)}}{\pi_1^{(3)}}\right)^{n_1}\left(\frac{\epsilon_1^{(4)}}{\epsilon_1^{(3)}}\right)^{a_1}-1$$

The numbers (complex or 17-adic) appearing in $\lambda$, besides $-\imath$ are now

$$\delta := \frac{\pi_1^{(4)}}{\pi_1^{(3)}} = \begin{cases} \frac{1}{17}(-109+11\vartheta^2-(37\vartheta-\frac{7}{2}\vartheta^3)\imath) & \text{in case of (42)} \\ \frac{1}{289}(1633-200\vartheta^2+(536\vartheta-66\vartheta^3)\imath) & \text{in case of (43)} \end{cases}$$

$$\beta := \frac{\epsilon^{(4)}}{\epsilon^{(3)}} = 12097-1320\vartheta^2+(3996\vartheta-436\vartheta^3)\imath$$

Now $\lambda = -\imath\,\delta^{n_1}\beta^{a_1}-1$ and, in analogy with section 6.2, we consider the 17-adic linear form $\Lambda_1 = \log_{17}(\lambda+1)$ and the complex linear form $\Lambda_0 = \imath^{-1}\mathrm{Log}(1+\lambda)$. On "expanding" Log in the right-hand side a new integer unknown $a_0$ makes its appearance as a coefficient of $2\pi = 2\mathrm{Log}(-1)$. Since $\imath^2=-1$, we have now $\imath\Lambda_0 = n_1\mathrm{Log}\,\delta + a_1\mathrm{Log}\,\beta + (1-4a_0)\mathrm{Log}(-\imath)$.

In complete analogy with the strategy of section 6.2, we next find in terms of $\log H = \log\max\{|a_1|,|a_2|,n_1\}$ an upper bound for $\mathrm{ord}_{17}(\lambda)$, by applying Yu's theorem and a lower bound for $|\Lambda_0|$ by applying the theorem of Baker and Wüstholz. From these we find upper bounds for $H$ and $n_1$ of very similar size with those found in section 6.2. An analogous reduction process to that described in section 6.3 has as a result the drastic reduction of $H$ down to 20 and of $n_1$ down to 5. A final check, similar to that in the end of section 6.3 shows that no solutions come from (42) and only the solution $(|u|,|v|,n_1)=(31,10,2)$ from (43). Since $z=2n_1$, this completes the proof of Proposition 7.1.

For the needs of the proof of Theorem 4.2 the $p$-adic first reduction step is sufficient, for, as a result of it an upper bound for $n_1$ less than 100 is obtained and then one has just to check whether $U_n(\pm4,-17)$ has the required shape just for those values of $n$ that are powers of 2 less than 200.

**Acknowledgement:** we are grateful to the anonymous referee for helpful suggestions to improve this paper.

Department of Mathematics,
Arizona State University,
Tempe AZ 85287-1804,
USA.
E-mail: bremner@asu.edu,
URL: math.asu.edu/˜andrew/bremner.html

Department of Mathematics,
University of Crete,
Iraklion,
Greece.
E-mail:tzanakis@math.uoc.gr
URL: www.math.uoc.gr/˜tzanakis